\documentclass[11pt,oneside]{amsart}

\usepackage{newtxtext}
\usepackage{newtxmath}
\usepackage{fix-cm}

\usepackage{tikz-cd}
\usepackage[arrow,curve,matrix,tips,2cell]{xy}

\usepackage[hmargin=2cm,vmargin=1.5cm]{geometry}

\usepackage{fancyhdr}

\usepackage{amsmath}
\usepackage{amscd}
\usepackage{amssymb}
\usepackage{latexsym}
\usepackage{url}

\usepackage{graphicx}

\newtheorem{theorem}{Theorem}[section]
\newtheorem*{theorem*}{Theorem}
\newtheorem{lemma}[theorem]{Lemma}
\newtheorem*{lemma*}{Lemma}
\newtheorem{corollary}[theorem]{Corollary}
\newtheorem{proposition}[theorem]{Proposition}

\newtheorem{remark}[theorem]{Remark}
\newtheorem{definition}[theorem]{Definition}
\newtheorem*{definition*}{Definition}
\newtheorem{question}[theorem]{Question}
\newtheorem*{question*}{Question}

\newtheorem{example}[theorem]{Example}
\newtheorem{examples}[theorem]{Examples}

\makeatletter
\def\revddots{\mathinner{\mkern1mu\raise\p@
\vbox{\kern7\p@\hbox{.}}\mkern2mu
\raise4\p@\hbox{.}\mkern2mu\raise7\p@\hbox{.}\mkern1mu}}
\makeatother % ddots gespiegelt
\newcommand{\bgl}{\begin{equation}} %eine Gleichung mit Ziffer
\newcommand{\egl}{\end{equation}}
\newcommand{\bgloz}{\begin{equation*}} %eine Gleichung ohne Ziffer
\newcommand{\egloz}{\end{equation*}}
\newcommand{\bgln}{\begin{eqnarray}} %mehrere Gleichungen mit Ziffer
\newcommand{\egln}{\end{eqnarray}}
\newcommand{\bglnoz}{\begin{eqnarray*}} %mehrere Gleichungen ohne Ziffer
\newcommand{\eglnoz}{\end{eqnarray*}}
\newcommand{\btheo}{\begin{theorem}}
\newcommand{\etheo}{\end{theorem}}
\newcommand{\btheooz}{\begin{theorem*}}
\newcommand{\etheooz}{\end{theorem*}}

\newcommand{\blemma}{\begin{lemma}}
\newcommand{\elemma}{\end{lemma}}
\newcommand{\blemmaoz}{\begin{lemma*}}
\newcommand{\elemmaoz}{\end{lemma*}}
\newcommand{\bproof}{\begin{proof}}
\newcommand{\eproof}{\end{proof}}
\newcommand{\bbew}{\begin{beweis}}
\newcommand{\ebew}{\end{beweis}}
\newcommand{\bremark}{\begin{remark}\em}
\newcommand{\eremark}{\end{remark}}
\newcommand{\bdefin}{\begin{definition}}
\newcommand{\edefin}{\end{definition}}
\newcommand{\bdefinoz}{\begin{definition*}}
\newcommand{\edefinoz}{\end{definition*}}
\newcommand{\bex}{\begin{example}}
\newcommand{\eex}{\end{example}}
\newcommand{\bexs}{\begin{examples}}
\newcommand{\eexs}{\end{examples}}

\newcommand{\bprop}{\begin{proposition}}
\newcommand{\eprop}{\end{proposition}}
\newcommand{\bcor}{\begin{corollary}}
\newcommand{\ecor}{\end{corollary}}
\newcommand{\bfa}{\begin{cases}} %Fallunterscheidung
\newcommand{\efa}{\end{cases}}

\newcommand{\bquestion}{\begin{question}}
\newcommand{\equestion}{\end{question}}
\newcommand{\bquestionoz}{\begin{question*}}
\newcommand{\equestionoz}{\end{question*}}
\newcommand{\cC}{\mathcal C}

\def\Cz{\mathbb{C}}

\def\Zz{\mathbb{Z}}

\def\1z{\mathbb{1}}
\newcommand{\ma}{\mapsto} % wird abgebildet auf
\newcommand{\onto}{\twoheadrightarrow} % surjektiv
\newcommand{\into}{\hookrightarrow} % injektiv
\def\SEMI{\mbox{$\times\kern-2pt\vrule height5pt width.6pt \kern3pt $}}

\newcommand{\reg}{^\times} % regulaer
\newcommand{\defeq}{\mathrel{:=}} % per Definition
\newcommand{\dop}{\text{: }} % in Mengen
\newcommand{\ilim}{\varinjlim} % induktiver Limes
\newcommand{\lge}{\left\{} % links geschweift
\newcommand{\rge}{\right\}} % rechts geschweift
\newcommand{\lru}{\left(} % links rund
\newcommand{\rru}{\right)} % rechts rund
\newcommand{\rukl}[1]{\lru #1 \rru} % runde Klammer
\newcommand{\gekl}[1]{\lge #1 \rge} % geschweifte Klammer
\newcommand{\menge}[2]{\gekl{ #1 \dop #2 }} % Menge
\begin{document}

\title{K-theory for generalized Lamplighter groups}

\thispagestyle{fancy}

\begin{abstract}
We compute K-theory for the reduced group C*-algebras of generalized Lamplighter groups.
\end{abstract}

\author{Xin Li}

\address{Xin Li, School of Mathematical Sciences, Queen Mary University of London, Mile End Road, London E1 4NS, United Kingdom}
\email{xin.li@qmul.ac.uk}


\maketitle

\setlength{\parindent}{0cm} \setlength{\parskip}{0.5cm}

\section{Introduction}

The classical Lamplighter group is given by the semidirect product $\rukl{\bigoplus_{\Zz} (\Zz / 2 \Zz)} \rtimes \Zz$,
where the $\Zz$-action on $\bigoplus_{\Zz} (\Zz / 2 \Zz)$ is induced by the canonical translation action of $\Zz$ on itself. This construction can be generalized by replacing $\Zz / 2 \Zz$ and $\Zz$ by other groups. The classical Lamplighter group and its generalizations are important examples in group theory which led to solutions of several open problems (see for instance \cite{Dyu,Dym,Gra}).

The goal of these notes is to derive a $K$-theory formula for group C*-algebras of generalized Lamplighter groups of the form $\rukl{\bigoplus_{\Gamma} \Sigma} \rtimes \Gamma$, where $\Sigma$ is an arbitrary finite group and $\Gamma$ is an arbitrary countable group. As in the classical setting, the $\Gamma$-action on $\bigoplus_{\Gamma} \Sigma$ is induced by the canonical left translation action of $\Gamma$ on itself. Our computations are inspired by \cite{FPV,Poo}, which treat the special case of free groups $\Gamma$ (\cite{FPV} deals with the case $\Gamma = \Zz$). Our method, however, is completely different from the ones adopted in \cite{FPV,Poo}.

Our main result reads as follows. Let $\Sigma$ be a finite group and $\Gamma$ a countable group. Let ${\rm con} \, \Sigma$ be the set of conjugacy classes in $\Sigma$, and ${\rm con}\reg \, \Sigma \defeq {\rm con} \, \Sigma \setminus \gekl{\gekl{1}}$ the set of non-trivial conjugacy classes. Let $\cC$ be the set of conjugacy classes of finite subgroups of $\Gamma$. For a finite subgroup $C$ of $\Gamma$, let $F(C)$ be the set of non-empty finite subsets of the right coset space $C \backslash \Gamma$ which are not of the form $\pi^{-1}(Y)$ for a finite subgroup $D \subseteq \Gamma$ with $C \subsetneq D$ and $Y \subseteq D \backslash \Gamma$, where $\pi: \: C \backslash \Gamma \onto D \backslash \Gamma$ is the canonical projection. The normalizer $N_C \defeq \menge{\gamma \in \Gamma}{\gamma C \gamma^{-1} = C}$ acts on $F(C)$ by left multiplication, and we denote the set of orbits by $N_C \backslash F(C)$. Given $X \in F(C)$, we write $C \cdot X \defeq \bigsqcup_{x \in X} C \cdot x$ and let $({\rm con}\reg \, \Sigma)^{C \cdot X}$ be the set of functions $C \cdot X \to {\rm con}\reg \, \Sigma$. $\gamma \in C$ acts on $\varphi \in ({\rm con}\reg \, \Sigma)^{C \cdot X}$ via $(\gamma.\varphi)(x) = \varphi(\gamma^{-1}x)$, and we set ${\rm Stab}_C(\varphi) = \menge{\gamma \in C}{\gamma.\varphi = \varphi}$ for $\varphi \in ({\rm con}\reg \, \Sigma)^{C \cdot X}$.
\btheo
\label{THM:main}
If $\Gamma$ satisfies the Baum-Connes conjecture with coefficients, then the K-theory of $C^*_{\lambda} ((\bigoplus_{\Gamma} \Sigma) \rtimes \Gamma)$ is given by
$$
  K_* (C^*_{\lambda} ((\bigoplus_{\Gamma} \Sigma) \rtimes \Gamma)) \cong K_*(C^*_{\lambda}(\Gamma)) \oplus \rukl{\bigoplus_{[C] \in \cC} \ \ \bigoplus_{[X] \in N_C \backslash F(C)} \ \ \bigoplus_{[\varphi] \in C \backslash \rukl{({\rm con}\reg \, \Sigma)^{C \cdot X}}} K_*(C^*_{\lambda}({\rm Stab}_C(\varphi)))}.
$$
Here we take one representative $C$ out of each class in $\cC$, one representative $X$ out of each class in $N_C \backslash F(C)$, and one representative $\varphi$ out of each class in $C \backslash \rukl{({\rm con}\reg \, \Sigma)^{C \cdot X}}$.
\etheo
\setlength{\parindent}{0cm} \setlength{\parskip}{0cm}

We refer the reader to \cite{BCH,Val,CELY} and the references therein for more information about the Baum-Connes conjecture. For instance, Theorem~\ref{THM:main} applies to all groups with the Haagerup property \cite{HK} and all hyperbolic groups \cite{Laf}.
\setlength{\parindent}{0cm} \setlength{\parskip}{0.5cm}

Note that $\Sigma$ enters our formula only in the form of ${\rm con}\reg \, \Sigma$. What is more, if $\Gamma$ is infinite, then for each $[C] \in \cC$, we simply get a free abelian group of  countably infinite rank, so that $K_* (C^*_{\lambda} ((\bigoplus_{\Gamma} \Sigma) \rtimes \Gamma))$ does not depend on $\Sigma$ at all. This becomes particularly evident in $K_1$, where Theorem~\ref{THM:main} yields the following
\bcor
\label{Cor:K1}
Let $\Sigma$ be a finite group and $\Gamma$ a countable group. If $\Gamma$ satisfies the Baum-Connes conjecture with coefficients, then the canonical inclusion $\Gamma \into \Sigma \rtimes \Gamma$ induces an isomorphism
$$
  K_1(C^*_{\lambda}(\Gamma)) \cong K_1(C^*_{\lambda}((\bigoplus_{\Gamma} \Sigma) \rtimes \Gamma)).
$$
\ecor

Moreover, for torsion-free $\Gamma$, our formula becomes particularly simple.
\bcor
\label{COR:torsionfree}
Let $\Sigma$ and $\Gamma$ be as in Theorem~\ref{THM:main}. Assume that $\Gamma$ is torsion-free. Write ${\rm FIN}\reg$ for the set of non-empty finite subsets of $\Gamma$. Then, under the same assumptions as in Theorem~\ref{THM:main}, we have
$$
  K_* (C^*_{\lambda} ((\bigoplus_{\Gamma} \Sigma) \rtimes \Gamma)) \cong K_*(C^*_{\lambda}(\Gamma)) \oplus \rukl{\bigoplus_{[X] \in \Gamma \backslash {\rm FIN}\reg} \ \ \bigoplus_{({\rm con}\reg \, \Sigma)^X} K_*(\Cz)}.
$$
\ecor

The proof of our main theorem proceeds in two steps. First, using the Going-Down principle from \cite{CEO,ENO} (see also \cite[\S~3]{CELY}), we show that $C^*_{\lambda} ((\bigoplus_{\Gamma} \Sigma) \rtimes \Gamma)$ has the same K-theory as the crossed product $C( ({\rm con} \, \Sigma)^{\Gamma} ) \rtimes_r \Gamma$ for the topological full shift $\Gamma \curvearrowright ({\rm con} \, \Sigma)^{\Gamma}$. Here we view ${\rm con} \, \Sigma$ as a finite alphabet. Secondly, we compute K-theory for $C( ({\rm con} \, \Sigma)^{\Gamma} ) \rtimes_r \Gamma$ using \cite{CEL1,CEL2}. As a by-product, we obtain a general K-theory formula for crossed products of topological full shifts (see Proposition~\ref{PROP:KFullShift}). Both steps require our assumption that $\Gamma$ satisfies the Baum-Connes conjecture with coefficients.

We point out that it is not possible to apply the results in \cite{CEL1,CEL2} directly because \cite{CEL1,CEL2} only deal with crossed products attached to actions on commutative C*-algebras.

I am indebted to Alain Valette for inviting me to Neuchatel, and to Sanaz Pooya and Alain Valette for interesting discussions which led to these notes. Moreover, I thank the referee for helpful comments and careful proofreading.

\section{K-theory for certain crossed products and generalized Lamplighter groups}

We first discuss the following abstract situation: Let $A = \bigoplus_{i=0}^n M_{k_i}$ be a finite dimensional C*-algebra, where $M_k$ is the algebra of $k \times k$-matrices. We assume that $k_0 = 1$, i.e., $A = \Cz \oplus M_{k_1} \oplus \dotso \oplus M_{k_n}$. Let $\Gamma$ be a countable group. We form the tensor product $\bigotimes_{\Gamma} A$ as follows: For every finite subset $F \subseteq \Gamma$, we form the ordinary tensor product $\bigotimes_F A$, and for $F_1 \subseteq F_2$, we have the canonical embedding $\bigotimes_{F_1} A \into \bigotimes_{F_2} A, \, x \ma x \otimes 1$ (here $1$ denotes the unit of $\bigotimes_{F_2 \setminus F_1} A$, and we used the canonical isomorphism $\bigotimes_{F_2} A \cong \rukl{\bigotimes_{F_1} A} \otimes \rukl{\bigotimes_{F_2 \setminus F_1} A}$). Then set $\bigotimes_{\Gamma} A \defeq \ilim_F \bigotimes_F A$. The left $\Gamma$-action on itself by translations induces an action $\Gamma \curvearrowright \bigotimes_{\Gamma} A$. Our goal is to compute the $K$-theory of $\rukl{\bigotimes_{\Gamma} A} \rtimes_r \Gamma$. The special case $A = C^*_{\lambda}(\Sigma)$ will lead to Theorem~\ref{THM:main}.

Let $e_i$ be a minimal projection in $M_{k_i} \subseteq A$. In particular, $e_0 = 1 \in \Cz \subseteq A$. For $F \subseteq \Gamma$ finite, let $\varphi \in \gekl{1, \dotsc, n}^F$, i.e., $\varphi$ is a function $\varphi: \: F \to \gekl{1, \dotsc, n}$. Define $e_{\varphi} \defeq \bigotimes_{f \in F} e_{\varphi(f)} \in \bigotimes_F A \subseteq \bigotimes_{\Gamma} A$. If $F = \emptyset$, then for $\varphi: \: \emptyset \to \gekl{1, \dotsc, n}$, we set $e_{\varphi} \defeq 1$ (where $1$ denotes the unit of $\bigotimes_{\Gamma} A$). The set
\begin{equation}
\label{e}
  \menge{e_{\varphi}}{\varphi \in \gekl{1, \dotsc, n}^F, \, F \subseteq \Gamma \ {\rm finite}}
\end{equation}
is a $\Gamma$-invariant family of commuting non-zero projections, which is closed under multiplication up to zero (i.e., the product of two projections in the family is either zero or again a projection in the family). We do not need it now, but the family is also linearly independent (see Lemma~\ref{LEM:D=FullShift} and the proof of \eqref{EQ:KFullShift--FIN}). Let $D$ be the C*-subalgebra of $\bigotimes_{\Gamma} A$ generated by the projections in \eqref{e}. Let $\iota: \: D \into \bigotimes_{\Gamma} A$ be the canonical embedding. Note that $\iota$ is $\Gamma$-equivariant.

\bprop
\label{PROP:D=A}
If $\Gamma$ satisfies the Baum-Connes conjecture with coefficients, then $\iota \rtimes_r \Gamma$ induces an isomorphism $K_*(D \rtimes_r \Gamma) \cong K_*(\rukl{\bigotimes_{\Gamma} A} \rtimes_r \Gamma)$.
\eprop
\setlength{\parindent}{0cm} \setlength{\parskip}{0cm}

\bproof
By the Going-Down principle (see \cite[\S~3]{CELY}), it suffices to show that for every finite subgroup $H \subseteq \Gamma$, $\iota \rtimes_r H$ induces an isomorphism $K_*(D \rtimes_r H) \cong K_*(\rukl{\bigotimes_{\Gamma} A} \rtimes_r H)$.
\setlength{\parindent}{0cm} \setlength{\parskip}{0.5cm}

Let us first treat the case of the trivial subgroup, $H = \gekl{1}$. For a fixed finite subset $F \subseteq \Gamma$, let
$$
D_F = C^*(\menge{e_{\varphi}}{\varphi \in \gekl{1, \dotsc, n}^{F'} \ {\rm for} \ F' \subseteq F}).
$$ 
Then $D = \ilim_F D_F$. We also have $\bigotimes_{\Gamma} A = \ilim_F \bigotimes_F A$. As $K$-theory is continuous, i.e., preserves direct limits, it suffices to show that $\iota_F \defeq \iota \vert_{D_F}: \: D_F \to \bigotimes_F A$ induces an isomorphism in $K_*$. Let $[\iota_F] \in KK(D_F,\bigotimes_F A)$ be the $KK$-element determined by $\iota_F$. Consider the projection $e = \sum_{i=0}^n e_i$ in $A$. $e$ is a full projection in $A$, and we have $eAe = \bigoplus_{i=0}^n \Cz e_i$. The $\bigotimes_F A$--$\bigotimes_F eAe$-imprimitivity bimodule $\bigotimes_F Ae$ gives rise to a $KK$-element $\mathbf{j}_F \in KK(\bigotimes_F A,\bigotimes_F eAe)$. $\mathbf{j}_F$ is invertible, and its inverse is the $KK$-element induced by the inclusion $\bigotimes_F eAe \into \bigotimes_F A$. Hence it suffices to show that the Kasparov product $[\iota_F] \cdot \mathbf{j}_F \in KK(D_F,\bigotimes_F eAe)$ induces an isomorphism $K_*(D_F) \to K_*(\bigotimes_F eAe)$.

First, consider the special case of a single element subset, $F = \gekl{f}$ for some $f \in \Gamma$. Let us write $D_f \defeq D_{\gekl{f}}$, $\iota_f \defeq \iota_{\gekl{f}}$ and $\mathbf{j}_f \defeq \mathbf{j}_{\gekl{f}}$. Since $D_f = \Cz \cdot 1 + \Cz e_1 + \dotso + \Cz e_n$ (where $1$ denotes the unit of $\bigotimes_{\Gamma} A$) and $eAe = \Cz e_0 \oplus \Cz e_1 \oplus \dotso \oplus \Cz e_n$, we can describe the map $K_*(D_F) \to K_*(\bigotimes_F eAe)$ induced by $[\iota_f] \cdot \mathbf{j}_f$ by the commutative diagram
\begin{equation*}
\begin{tikzcd}
K_*(D_f) \ar{r} \ar[equal]{d} & K_*(eAe) \ar[equal]{d} \\
\Zz [1] \oplus \bigoplus_{i=1}^n \Zz [e_i] \ar["M_f"]{r} & \bigoplus_{i=0}^n \Zz [e_i]
\end{tikzcd}
\end{equation*}
where the upper horizontal map is the map we want to describe, and $M_f$ is the $(n+1) \times (n+1)$-matrix
\begin{equation*}
  M_f = 
  \begin{pmatrix}
  1 & 0 & \dotso & 0 \\
  k_1 & 1 & & 0 \\
  \vdots & & \ddots & \\
  k_n & 0 & & 1
  \end{pmatrix}.
\end{equation*}
Obviously, $M_f$ is invertible. Note that everything is independent of $f$.

Now consider the case of a general finite subset $F \subseteq \Gamma$. Since $D_F = \bigotimes_{f \in F} D_f$, we have $K_*(D_F) \cong \bigotimes_{f \in F} K_*(D_f)$, and we also have $K_*(\bigotimes_F eAe) \cong \bigotimes_{f \in F} K_*(eAe)$. The homomorphism $K_*(D_F) \to K_*(\bigotimes_F eAe)$ induced by
$[\iota_F] \cdot \mathbf{j}_F$ respects this tensor product decomposition, in the sense that we have a commutative diagram
\begin{equation*}
\begin{tikzcd}
\bigotimes_{f \in F} K_*(D_f) \ar[equal]{d} \ar[equal,"\sim"]{r} & K_*(D_F) \ar{r} & K_*(\bigotimes_F eAe) \ar[equal,"\sim"]{r} & \bigotimes_{f \in F} K_*(eAe) \ar[equal]{d} \\
\bigotimes_{f \in F} \rukl{\Zz [1] \oplus \bigoplus_{i=1}^n \Zz [e_i]} \ar["M_F \, = \, \otimes_{f \in F} M_f"]{rrr} & & & \bigotimes_{f \in F} \rukl{\bigoplus_{i=0}^n \Zz [e_i]}
\end{tikzcd}
\end{equation*}
Again, we see that $M_F$ is invertible because all the $M_f$, $f \in F$, are.

Now let us deal with the case of an arbitrary finite subgroup $H \subseteq \Gamma$. If we choose an increasing sequence of $H$-invariant finite subsets $F \subseteq \Gamma$ whose union is $\Gamma$, we obtain $H$-equivariant inductive limit decompositions $D = \ilim_F D_F$ and $\bigotimes_{\Gamma} A = \ilim_F \bigotimes_F A$. Hence, again by continuity of $K$-theory, it suffices to show that, for every $F$, $\iota_F \rtimes_r H: \: D_F \rtimes_r H \to \rukl{\bigotimes_F A} \rtimes_r H$ induces an isomorphism in $K_*$. Let $\mathbf{j}_F \in KK(\bigotimes_F A,\bigotimes_F eAe)$ be as before. Since the full projection $\bigotimes_F e \in \bigotimes_F A$ giving rise to $\mathbf{j}_F$ is $H$-invariant, $\mathbf{j}_F$ is a $KK^H$-equivalence (see \cite[Remark~3.3.16]{CELY}). Thus, to show that $\iota_F \rtimes_r H: \: D_F \rtimes_r H \to \rukl{\bigotimes_F A} \rtimes_r H$ induces an isomorphism in $K_*$, it suffices to show that $[\iota_F] \cdot \mathbf{j}_F \in KK^H(D_F,\bigotimes_F eAe)$ induces an isomorphism $K_*(D_F \rtimes_r H) \to K_*(\rukl{\bigotimes_F eAe} \rtimes_r H)$, for which in turn it is enough to prove that $[\iota_F] \cdot \mathbf{j}_F$ is a $KK^H$-equivalence.

Now both $D_F$ and $\bigotimes_F eAe$ are finite dimensional commutative C*-algebras with an $H$-action, so that we are exactly in the setting of \cite[Appendix]{CEL2}. It is straightforward to check that $[\iota_F] \cdot \mathbf{j}_F = x^H_{M_F}$, where $x^H_{M_F}$ is the element in $KK^H(D_F,\bigotimes_F eAe)$ corresponding to the matrix $M_F$, as constructed in \cite[Appendix]{CEL2}. By \cite[Lemma~A.2]{CEL2}, $x^H_{M_F}$ is a $KK^H$-equivalence because $M_F$ is an invertible matrix. The inverse of $x^H_{M_F}$ is given by $x^H_{M_F^{-1}}$.
\eproof
\setlength{\parindent}{0cm} \setlength{\parskip}{0.5cm}

\bremark
Note that our assumption on $A$ that $\Cz$ appears as a direct summand is really necessary. For instance, if $A = M_2$, then $\bigotimes_{\Gamma} A$ would be the UHF algebra $M_{2^{\infty}}$ (as soon as $\Gamma$ is infinite). But we have $K_0(M_{2^{\infty}}) \cong \Zz[\frac{1}{2}]$, while our method would always yield a free abelian group for $K_0$. Hence our method fails.
\eremark

Let us now compare with the topological full shift $\Gamma \curvearrowright \gekl{0, \dotsc, n}^{\Gamma}$. For a finite subset $F \subseteq \Gamma$, let $\pi_F$ be the canonical projection $\gekl{0, \dotsc, n}^{\Gamma} \onto \gekl{0, \dotsc, n}^F$. Given $\varphi \in \gekl{0, \dotsc, n}^F$, we have the cylinder set $\pi_F^{-1}(\varphi)$ and its characteristic function $1_{\pi_F^{-1}(\varphi)} \in C(\gekl{0, \dotsc, n}^{\Gamma})$. The following is now easy to see.
\blemma
\label{LEM:D=FullShift}
The $\Gamma$-equivariant isomorphism $D \cong C(\gekl{0, \dotsc, n}^{\Gamma}), \, e_{\varphi} \mapsto 1_{\pi_F^{-1}(\varphi)}$ induces an isomorphism $D \rtimes_r \Gamma \cong C(\gekl{0, \dotsc, n}^{\Gamma}) \rtimes_r \Gamma$.
\elemma

We now compute K-theory for $C(\gekl{0, \dotsc, n}^{\Gamma}) \rtimes_r \Gamma$.
\bprop
\label{PROP:KFullShift}
If $\Gamma$ satisfies the Baum-Connes conjecture with coefficients, then
$$
  K_*(C(\gekl{0, \dotsc, n}^{\Gamma}) \rtimes_r \Gamma) \cong K_*(C^*_{\lambda}(\Gamma)) \oplus \rukl{\bigoplus_{[C] \in \cC} \ \ \bigoplus_{[X] \in N_C \backslash F(C)} \ \ \bigoplus_{[\varphi] \in C \backslash \rukl{\gekl{1, \dotsc, n}^{C \cdot X}}} K_*(C^*_{\lambda}({\rm Stab}_C(\varphi)))}.
$$
Here we use the same notation as in Theorem~\ref{THM:main}.
\eprop
\setlength{\parindent}{0cm} \setlength{\parskip}{0cm}

\bproof
First of all, the same arguments as for \cite[Examples~2.13 \& 3.1]{CEL2} show that the family
\begin{equation*}
  \menge{\pi_F^{-1}(\varphi)}{\varphi \in \gekl{1, \dotsc, n}^F, \, F \subseteq \Gamma \ {\rm finite}}
\end{equation*}
is a $\Gamma$-invariant regular basis for the compact open sets in $\gekl{0, \dotsc, n}^{\Gamma}$. Here $\Gamma$ acts via $\gamma.\pi_F^{-1}(\varphi) = \pi_{\gamma \cdot F}^{-1}(\gamma.\varphi)$, where $\gamma.\varphi \in \gekl{1, \dotsc, n}^{\gamma \cdot F}$ is given by $(\gamma.\varphi)(x) = \varphi(\gamma^{-1}x)$. Therefore, using the bijection
$$
  \bigsqcup_{[F] \in \Gamma \backslash {\rm FIN}} {\rm Stab}_{\Gamma}(F) \backslash \rukl{\gekl{1, \dotsc, n}^F} \cong \Gamma \, \backslash \menge{\pi_F^{-1}(\varphi)}{\varphi \in \gekl{1, \dotsc, n}^F, \, F \subseteq \Gamma \ {\rm finite}}, \, [\varphi] \ma [\varphi],
$$ 
and the observation that for $\gamma \in \Gamma$ and $\varphi \in \gekl{1, \dotsc, n}^F$, we have $\gamma.\pi_F^{-1}(\varphi) = \pi_{F}^{-1}(\varphi)$ if and only if $\gamma \cdot F = F$ and $\gamma.\varphi = \varphi$, we may apply \cite[Corollary~3.14]{CEL2}, and obtain
\begin{equation}
\label{EQ:KFullShift--FIN}
  K_*\rukl{C\rukl{\gekl{0, \dotsc, n}^{\Gamma}} \rtimes_r \Gamma} \cong \bigoplus_{[F] \in \Gamma \backslash {\rm FIN}} \ \ \bigoplus_{[\varphi] \in {\rm Stab}_{\Gamma}(F) \backslash \rukl{\gekl{1, \dotsc, n}^F}} K_*(C^*_{\lambda}({\rm Stab}_{\Gamma}(\varphi))).
\end{equation}
Here ${\rm FIN}$ is the set of all finite subsets of $\Gamma$, and $\Gamma \backslash {\rm FIN}$ is the set of orbits of the left translation action $\Gamma \curvearrowright {\rm FIN}$. Moreover, ${\rm Stab}_{\Gamma}(F)$ and ${\rm Stab}_{\Gamma}(\varphi)$ denote the stabilizer groups ${\rm Stab}_{\Gamma}(F) = \menge{\gamma \in \Gamma}{\gamma \cdot F = F}$ and ${\rm Stab}_{\Gamma}(\varphi) = \menge{\gamma \in \Gamma}{\gamma.\varphi = \varphi}$, and $C^*_{\lambda}({\rm Stab}_{\Gamma}(\varphi))$ is the reduced group C*-algebra of ${\rm Stab}_{\Gamma}(\varphi)$.
\setlength{\parindent}{0cm} \setlength{\parskip}{0.5cm}

Now we analyse $\Gamma \backslash {\rm FIN}$ and ${\rm Stab}_{\Gamma}(F)$ for $[F] \in \Gamma \backslash {\rm FIN}$. For $F = \emptyset$, we have ${\rm Stab}_{\Gamma}(\varphi) = \Gamma$. This yields $K_*(C^*_{\lambda}(\Gamma))$ as one direct summand on the right-hand side of \eqref{EQ:KFullShift--FIN}. We set ${\rm FIN}\reg \defeq {\rm FIN} \setminus \gekl{\emptyset}$. Let us describe $\Gamma \backslash {\rm FIN}\reg$. Let $\cC$, $F(C)$ and $N_C$ be as in Theorem~\ref{THM:main}. Then we claim that
\begin{equation}
\label{EQ:NCF=}
  \bigsqcup_{[C] \in \cC} N_C \backslash F(C) \to \Gamma \backslash {\rm FIN}\reg, \, [X] \ma [C \cdot X]
\end{equation}
is a bijection, and that for every $[C] \in \cC$, $[X] \in N_C \backslash F(C)$, we have
\begin{equation}
\label{EQ:Stab=C}
  {\rm Stab}_{\Gamma}(C \cdot X) = C.
\end{equation}
First note that the map \eqref{EQ:NCF=} is well-defined. Moreover, this map is surjective because every $F \in {\rm FIN}\reg$ with ${\rm Stab}_{\Gamma}(F) = C$ is of the form $F = C \cdot X$ for some finite, non-empty subset $X \subseteq C \backslash \Gamma$. Now, $X$ must lie in $F(C)$. Suppose not, i.e., $X = \pi^{-1}(Y)$ for a finite subgroup $D \subseteq \Gamma$ with $C \subsetneq D$ and $Y \subseteq D \backslash \Gamma$, where $\pi: \: C \backslash \Gamma \onto D \backslash \Gamma$ is the canonical projection. Then $F = C \cdot X = D \cdot Y$, so that $D \subseteq {\rm Stab}_{\Gamma}(F)$, in contradiction to ${\rm Stab}_{\Gamma}(F) = C$. Not only does this show surjectivity, but it also proves \eqref{EQ:Stab=C}. To see injectivity of \eqref{EQ:NCF=}, assume that $X \in F(C)$ and $X' \in F(C')$ satisfy $[C \cdot X] = [C' \cdot X']$, say $C' \cdot X' = \gamma \cdot C \cdot X$. It follows that $C' = {\rm Stab}_{\Gamma}(C' \cdot X') = \gamma {\rm Stab}_{\Gamma}(C \cdot X) \gamma^{-1} = \gamma C \gamma^{-1}$. Hence $[C] = [C']$, and since we are taking one representative out of each class, we must actually have $C = C'$. But then $\gamma$ must lie in $N_C$, and we must have $C \cdot X' = \gamma \cdot C \cdot X = C \cdot \gamma \cdot X$, so that $X' = \gamma \cdot X$, i.e., $[X'] = [X]$ in $N_C \backslash F(C)$. This shows injectivity.

We now complete the proof by plugging the bijections \eqref{EQ:NCF=}, \eqref{EQ:Stab=C} into \eqref{EQ:KFullShift--FIN} and observing that for $X \in F(C)$ and $\varphi \in \gekl{1, \dotsc, n}^{C \cdot X}$, we have ${\rm Stab}_{\Gamma}(\varphi) \subseteq {\rm Stab}_{\Gamma}(C \cdot X) = C$.
\eproof
\setlength{\parindent}{0cm} \setlength{\parskip}{0.5cm}

Combining Proposition~\ref{PROP:D=A}, Lemma~\ref{LEM:D=FullShift} and Proposition~\ref{PROP:KFullShift}, and using the concrete construction in \cite[\S~3]{CEL2} for our following assertion on $K_1$, we obtain
\bcor
\label{COR:K(A}
In the situation of Proposition~\ref{PROP:D=A}, if $\Gamma$ satisfies the Baum-Connes conjecture with coefficients, then we have
$$
  K_*((\bigotimes_{\Gamma} A) \rtimes_r \Gamma) \cong K_*(C^*_{\lambda}(\Gamma)) \oplus \rukl{\bigoplus_{[C] \in \cC} \ \ \bigoplus_{[X] \in N_C \backslash F(C)} \ \ \bigoplus_{[\varphi] \in C \backslash \rukl{\gekl{1, \dotsc, n}^{C \cdot X}}} K_*(C^*_{\lambda}({\rm Stab}_C(\varphi)))}.
$$
In $K_1$, the canonical map $C^*_{\lambda}(\Gamma) \into (\bigotimes_{\Gamma} A) \rtimes_r \Gamma$ induces an isomorphism
$$
  K_1(C^*_{\lambda}(\Gamma)) \cong K_1((\bigotimes_{\Gamma} A) \rtimes_r \Gamma).
$$
If $\Gamma$ is in addition torsion-free, then we obtain
$$
  K_*((\bigotimes_{\Gamma} A) \rtimes_r \Gamma) \cong K_*(C^*_{\lambda}(\Gamma)) \oplus \rukl{\bigoplus_{[X] \in \Gamma \backslash {\rm FIN}\reg} \ \ \bigoplus_{\gekl{1, \dotsc, n}^X} K_*(\Cz)}.
$$
\ecor

Now let us apply our $K$-theory formula to generalized Lamplighter groups. Consider the case $A = C^*_{\lambda}(\Sigma)$ for a finite group $\Sigma$. Our assumption on $A$ that $\Cz$ appears as a direct summand is satisfied because the trivial representation gives rise to a direct summand $\Cz$ in $C^*_{\lambda}(\Sigma)$. The remaining direct summands of $A$ are in one-to-one correspondence with ${\rm con}\reg \, \Sigma$. Hence we obtain
\bcor
\label{COR:KFullShift}
Let $\Sigma$ be a finite group. If $\Gamma$ satisfies the Baum-Connes conjecture with coefficients, then we have 
$$
  K_* (C^*_{\lambda} ((\bigoplus_{\Gamma} \Sigma) \rtimes \Gamma)) \cong K_*(C^*_{\lambda}(\Gamma)) \oplus \rukl{\bigoplus_{[C] \in \cC} \ \ \bigoplus_{[X] \in N_C \backslash F(C)} \ \ \bigoplus_{[\varphi] \in C \backslash \rukl{({\rm con}\reg \, \Sigma)^{C \cdot X}}} K_*(C^*_{\lambda}({\rm Stab}_C(\varphi)))}.
$$
In $K_1$, the canonical inclusion $\Gamma \into \Sigma \rtimes \Gamma$ induces an isomorphism
$$
  K_1(C^*_{\lambda}(\Gamma)) \cong K_1(C^*_{\lambda}((\bigoplus_{\Gamma} \Sigma) \rtimes \Gamma)).
$$
If $\Gamma$ is in addition torsion-free, then we obtain
$$
  K_* (C^*_{\lambda} ((\bigoplus_{\Gamma} \Sigma) \rtimes \Gamma)) \cong K_*(C^*_{\lambda}(\Gamma)) \oplus \rukl{\bigoplus_{[X] \in \Gamma \backslash {\rm FIN}\reg} \ \ \bigoplus_{({\rm con}\reg \, \Sigma)^X} K_*(\Cz)}.
$$
\ecor
\setlength{\parindent}{0cm} \setlength{\parskip}{0cm}

This completes the proofs of Theorem~\ref{THM:main}, Corollary~\ref{Cor:K1} and Corollary~\ref{COR:torsionfree}.
\setlength{\parindent}{0cm} \setlength{\parskip}{0.5cm}

\bremark
As in \cite[Corollary~3.14]{CEL2}, we get $KK$-equivalences in Proposition~\ref{PROP:KFullShift}, Corollary~\ref{COR:K(A} and Corollary~\ref{COR:KFullShift} if $\Gamma$ satisfies the strong Baum-Connes conjecture. 
\eremark

\bremark
Moreover, as in \cite[Corollary~3.14]{CEL2}, we could allow coefficients in Proposition~\ref{PROP:KFullShift}, Corollary~\ref{COR:K(A} and Corollary~\ref{COR:KFullShift}. However, in Corollary~\ref{COR:KFullShift}, we would only get a $K$-theory formula for crossed products $B \rtimes_r \rukl{(\bigoplus_{\Gamma} \Sigma) \rtimes \Gamma}$ where the $(\bigoplus_{\Gamma} \Sigma)$-action on the C*-algebra $B$ is trivial.
\eremark

\end{document}